\documentclass{article}
\usepackage{arxiv}
\usepackage[utf8]{inputenc}
\usepackage[T2A]{fontenc}
\usepackage[english,russian]{babel}
\usepackage{mmap}
\usepackage{hyperref}
\usepackage{url}
\usepackage{mathtools}
\usepackage{amsfonts}
\usepackage{amsmath}
\usepackage{amssymb}
\usepackage{amsthm}
\usepackage{algorithm}
\usepackage{algorithmic}
\usepackage{graphicx}
\usepackage{titlesec}
\usepackage[labelsep=period]{caption}
\floatname{algorithm}{Алгоритм}
\hypersetup{
  colorlinks   = true,    
  urlcolor     = green,   
  linkcolor    = blue,    
  citecolor    = red      
}

\def\keywordname{{\bfseries Ключевые слова:}}%
\def\keywords#1{\par\addvspace\medskipamount{\rightskip=0pt plus1cm
\def\and{\ifhmode\unskip\nobreak\fi\ $\cdot$ }\noindent\keywordname\enspace\ignorespaces#1\par}}

\def\keywordnameeng{{\bfseries Keywords:}}%
\def\keywordseng#1{\par\addvspace\medskipamount{\rightskip=0pt plus1cm
\def\and{\ifhmode\unskip\nobreak\fi\ $\cdot$ }\noindent\keywordnameeng\enspace\ignorespaces#1\par}}

\titlelabel{\thetitle. }
\mathtoolsset{showonlyrefs}
\newcommand*{\No}{\textnumero}

\theoremstyle{plain}

\theoremstyle{definition}

\newtheorem*{task*}{Задача}

\renewcommand{\phi}{\varphi}
\renewcommand{\epsilon}{\varepsilon}

\renewcommand{\leq}{\leqslant}
\renewcommand{\geq}{\geqslant}

\DeclarePairedDelimiter{\rbr}{(}{)}

\DeclarePairedDelimiter{\cbr}{\{}{\}}
\DeclarePairedDelimiter{\abr}{\langle}{\rangle}
\DeclarePairedDelimiter{\abs}{\lvert}{\rvert}
\DeclarePairedDelimiter{\cardi}{\lvert}{\rvert}
\DeclarePairedDelimiter{\norm}{\lVert}{\rVert}

\newcommand*{\RR}{\mathbb{R}}
\renewcommand*{\emptyset}{\varnothing}

\newtheorem{teo}{Теорема}

\newtheorem{lem}{Лемма}
\newtheorem{cor}{Следствие}

\newtheorem{fed}{Определение}
\newtheorem{rem}{Замечание}

\newtheorem{exl}{Пример}

\renewcommand{\paragraph}{\section*}

\begin{document}

\title{Cубградиентные методы с шагом типа Б.\,Т.~Поляка для задач минимизации квазивыпуклых функций с ограничениями-неравенствами и аналогами острого минимума}

\author{
    Пучинин Сергей Максимович \\
    МФТИ \\
    Москва, Россия \\
    \texttt{puchinin.sm@phystech.edu} \\
    \And
    Корольков Егор Романович \\
    МФТИ \\
    Москва, Россия \\
    \texttt{korolkov.er@phystech.edu} \\
    \And
    Стонякин Федор Сергеевич \\
    МФТИ \\
    Москва, Россия \\
    КФУ им. В.\,И.\,Вернадского \\
    Симферополь, Россия \\
    \texttt{fedyor@mail.ru}
    \And
    Алкуса Мохаммад Соуд \\
    МФТИ \\
    Москва, Россия \\
    \texttt{mohammad.alkousa@phystech.edu} \\
    \And
    Выгузов Александр Альбертович \\
    МФТИ \\
    Москва, Россия \\
    \texttt{vyguzov.aa@phystech.edu} \\
}

{\raggedright УДК~519.85}\bigskip

\date{}

\maketitle

\begin{abstract}
\end{abstract}

\paragraph{1. Введение}

Задачи математического программирования возникают в самых разных приложениях. Одним из хорошо известных подходов к решению задач минимизации функций с ограничениями-неравенствами являются так называемые схемы с переключением между продуктивными и непродуктивными шагами, впервые предложенные Поляком в его работе \cite{Polyak1967}. Общая идея подхода заключается в следующем: если в текущей точке значение функции-ограничения достаточно хорошее, то спуск выполняется по целевой функции, а в противном случае~--- по функции ограничения. Такого типа подходам посвящаются всё новые работы как для выпуклых задач большой и сверхбольшой размерности \cite{Huang2023, Bayandina2018, Lagae2017, Nesterov2014}, так и для некоторых классов невыпуклых задач \cite{Huang2023}. 

В данной статье рассмотрены аналоги субградиентных методов такого типа с вариантом шага Б.\,Т.~Поляка c использованием информации о минимальном значении целевой функции $f^*$, хорошо известного для задач без ограничений \cite{Polyak1969}. Использование такого типа шага для липшицевых задач с острым минимумом позволяет доказать сходимость субградиентного метода со скоростью геометрической прогрессии. При этом важно, что реализация метода не потребует знания параметра острого минимума. Для других вариантов регулировки шага субградиентного метода такую скорость сходимости оказывается возможным доказать лишь за счёт процедуры рестартов по параметру острого минимума, оценка которого проблематична для многих реально возникающих типов прикладных задач. Стоит отметить, что использование в методах градиентного типа шага Б.\,Т.~Поляка довольно популярно (см., например, \cite{Hazan2019, Loizou2021}). Так, совсем недавно предложены и изучены различные современные варианты шага Поляка \cite{Wang2023, Devanathan2023, Abdukhakimov2023}.

В статье \cite{Ablaev2023} для липшицевых задач математического программирования предложены схемы рестартов по параметру <<условного>> острого минимума, для которых гарантирована сходимость в окрестность точного решения за линейное время как для выпуклых, так и для квазивыпуклых задач. Однако представляется существенным моментом требование знать значение параметра острого минимума для реализации метода. Цель предлагаемой статьи~--- обойти эту проблему посредством исследования субградиентных схем с переключениями по продуктивным и непродуктивным шагам с вариантом шага Б.\,Т.~Поляка. Однако наличие ограничений-неравенств приводит к нетривиальности вопроса оценки искомого точного минимального значения целевой функции. Более того, это значение может быть большим по сравнению с глобальным минимумом целевой функции на всём пространстве. Поэтому в настоящей статье исследуются аналоги шага Б.\,Т.~Поляка с неточной информацией о минимуме целевой функции и возможность гарантировать сходимость метода со скоростью геометрической прогрессии.

Статья состоит из введения, 4 основных разделов и заключения.

Во втором разделе приводятся основные определения и постановка рассматриваемой задачи, а также некоторые известные прикладные результаты, которые потребуются для анализа предлагаемых в последующих разделах алгоритмов.

В третьем разделе рассматривается аналог острого минимума, а именно $\epsilon$-острый минимум из работы \cite{Stonyakin2023}. Для данного варианта острого минимума предложен алгоритм, гарантирующий линейную скорость сходимости к $\epsilon$-точному решению задачи по функции и по ограничению.

В четвертом разделе рассматривается альтернативный вариант задания острого минимума, содержащий максимум из невязки по функции и значения ограничения \cite{Lin2020, Ablaev2023}. Данный подход к острому минимуму не предполагает задания параметра точности $\epsilon$. Это позволяет предложить алгоритм с гарантией линейной скорости сходимости в окрестность множества точных решений, размер которой определяется степенью точности информации о минимальном значении функции. Сходимость к точному решению со скоростью геометрической прогрессии  можно гарантировать только при условии доступности информации о точном значении минимума целевой функции, чего не удаётся доказать в случае использования $\epsilon$-острого минимума.

В пятом разделе приведены результаты численных экспериментов и продемонстрирована эффективность предлагаемых алгоритмов для некоторых конкретных задач. Среди них можно выделить задачу геометрического программирования \cite{boyd2007tutorial} и задачу проектирования механических конструкций \cite{nesterov_2014}.

\paragraph{2. Постановка задачи и основные теоретические сведения}

В работе рассматривается задача минимизации вида
\begin{gather}\label{eq:problem_statement}
    f(x) \to \min_{g(x) \leq 0}, \quad x \in Q,
\end{gather}
где $Q \subset \RR^n$~--- замкнутое выпуклое множество; $f: Q \longrightarrow \RR$~--- квазивыпуклая функция, удовлетворяющая условию Липшица для некоторой константы $M_f > 0 $ (то есть $\abs*{f(x) - f(y)} \leq M_f \norm*{x - y}_2$ для любых $x,y \in Q$); $g: Q \longrightarrow \RR$~--- $\mu$-слабо выпуклая функция, удовлетворяющая условию Липшица для некоторой константы $M_g > 0 $.

Напомним, что функция $g:\mathbb{R}^{n} \longrightarrow \mathbb{R}$ (аналогично и для функций $g: Q \longrightarrow \mathbb{R}$) называется $\mu$-слабо выпуклой $(\mu\geqslant 0)$, если функция $x \mapsto  g(x) +\frac{\mu}{2}\|x\|_2^{2}$ выпукла. Для недифференцируемых $\mu$-слабо выпуклых функций $g$ под субдифференциалом $\partial g(x)$ в точке $x$ можно понимать (см. \cite{Davis2018Subgradient} и цитируемую там литературу) множество всех векторов $\upsilon\in \mathbb{R}^{n}$, удовлетворяющих неравенству
\begin{equation}\label{f_1}
g(y)\geqslant g(x)+\langle\upsilon, y-x \rangle + o(\|y-x\|_2)\quad \text{при}\;\; y \rightarrow x.
\end{equation}
Известно \cite{Davis2018Subgradient}, что все векторы-субградиенты $\upsilon\in \mathbb{R}^{n}$ из \eqref{f_1} автоматически удовлетворяют неравенству
\begin{equation}\label{f_2}
g(y)\geqslant g(x)+\langle\upsilon, y-x\rangle  - \frac{\mu}{2}\|y-x\|_2^{2}\quad \forall\;x,y\in \mathbb{R}^{n}, \;\upsilon\in \partial g(x).
\end{equation}
Можно проверить, что слабо выпуклые функции локально липшицевы, и поэтому в качестве субградиентов можно использовать произвольный вектор из субдифференциала Кларка (см., например, \cite{Dudov2021}). 
Если целевая функция или функция ограничений недифференцируемы в некоторой точке $x$, то ввиду сделанного выше предположения об их липшицевости будем считать их субдифференцируемыми по Кларку и под $\nabla f(x)$ ($\nabla g(x)$) понимать произвольный ненулевой элемент субдифференциала Кларка в данной точке. Если же $f$($g$) дифференцируемы в обычном смысле, то под $\nabla f(x)$ ($\nabla g(x)$) будем понимать обычный градиент $f$ ($g$) в точке $x$.    

\begin{fed}[см.~{\cite[п.~1.5]{Nesterov1989}}]\label{fed:quasi1}
    Функция $f$, определенная на выпуклом множестве $Q \subset \RR^n$, называется квазивыпуклой на этом множестве, если
    \begin{gather}
        f\rbr*{\lambda x + (1 - \lambda) y} \leq \max\cbr*{f(x), f(y)} \quad \forall x, y \in Q \quad \forall \lambda \in [0; 1].
    \end{gather}
\end{fed}

Будем обозначать через $f^*$ оптимальное значение целевой функции задачи~\eqref{eq:problem_statement}, а через $X_*$~--- множество таких $x \in Q$, на которых это значение достигается. Также для произвольных $x, y \in Q$ введём вспомогательную величину
\begin{gather}\label{eq:v_f}
    v_f(x, y) \coloneqq \begin{cases}
        \abr*{\frac{\nabla f(x)}{\norm*{\nabla f(x)}_2}, x - y}, \quad &\nabla f(x) \neq 0;\\
        0, \quad &\nabla f(x) = 0,
    \end{cases}
\end{gather}
и докажем для неё следующую техническую лемму:
\begin{lem}\label{lem:v_f}
    Пусть $f: Q \longrightarrow \RR$ квазивыпукла и $M_f$-липшицева, тогда для любого $x \in Q$ такого, что $\nabla f(x) \neq 0$, выполнено
    \begin{gather}\label{eq:lem_v_f}
        f(x) - f^* \leq M_f v_f(x, x_*) \quad \forall x_* \in X_*.
    \end{gather}
\end{lem}
\proof* В силу \cite[сл.~1.5.1]{Nesterov1989}, а также условия $M_f$-липшицевости функции $f$ для произвольного $x_* \in X_*$
    \begin{gather}
        f(x) - f^* \leq \max\cbr*{f(y) - f(z) \mid \norm*{y - z}_2 \leq v_f(x, x_*)} \leq M_f v_f(x, x_*),
    \end{gather}
    что и требовалось.
\qed

Будем в дальнейшем предполагать, что у функции $f$ нет стационарных точек, за исключением, может быть, точек из $X_*$ (замечание 4, \cite{Stonyakin2023}). Легко видеть, что неравенство~\eqref{eq:lem_v_f} верно также и для точек $x$ из $X_*$, несмотря на возможное нулевое значение градиента. Тем самым дополнительное ограничение на $x$ ($\nabla f(x) \neq 0$) в условии леммы~\ref{lem:v_f} может быть опущено.

Зачастую на практике точное значение $f^*$ неизвестно. Поэтому будем считать, что алгоритмам, предложенным в данной работе, доступно некоторое значение $\overline{f}$, связанное с $f^*$ следующим образом:
\begin{gather}\label{eq:f_plus}
    f(x) - \overline{f} = c(x) \rbr*{f(x) - f^*},
\end{gather}
где $c(x) \in [C, 2 - C]$, $C \in (0; 1]$~--- некоторая константа. Очевидно, что равенство~\eqref{eq:f_plus} не может быть верно ни при каком значении $f(x)$, достаточно близком к $\overline{f}$ или $f^*$, где допустимая степень близости задается параметром $C$~--- чем он меньше, тем для более близких значений равенство~\eqref{eq:f_plus} выполняется. Все последующие исследования предлагаемых далее алгоритмов подразумевают, что значения функции в рассматриваемых последовательностях точек $x^k$ не приближаются к $\overline{f}$ или $f^*$ слишком близко, то есть для них равенство~\eqref{eq:f_plus} подразумевается верным. Однако при $C = 1$ равенство~\eqref{eq:f_plus}, очевидно, всегда верно, и никаких дополнительных ограничений не накладывается.

Также для анализа приведенных ниже алгоритмов нам потребуется следующее утверждение (см.~\cite{Polyak1969}):
\begin{lem}\label{lem:proj}
Пусть $\phi: Q \longrightarrow \RR$~--- некоторая функция ($Q \subset \RR^n$~--- замкнутое выпуклое множество), тогда для любых $x \in Q$, $x_* \in X_* \subset Q$ и любого $h \geq 0$ верно следующее неравенство:
\begin{gather}
    \norm*{\Bar{x} - x_*}_2^2 \leq \norm*{x - x_*}_2^2 - 2h \abr*{\nabla \phi(x), x - x_*} + h^2 \norm*{\nabla \phi(x)}_2^2,
\end{gather}
где $\Bar{x} = \operatorname{Pr}_Q \cbr*{x - h \nabla \phi(x)}$.
\end{lem}

\paragraph{3. Cубградиентный метод с шагом типа Б.\,Т.~Поляка для липшицевых задач с ограничениями-неравенствами с $\epsilon$-острым минимумом}

В данном разделе будет рассмотрен вариант острого минимума и cубградиентного метода, зависящий от желаемой точности решения $\varepsilon > 0$ \cite{Stonyakin2023}. Такой подход достаточно общий, но позволяет даже при доступности точной информации об $f^*$ гарантировать лишь попадание в окрестность некоторого точного решения задачи.

\begin{fed}
    Будем говорить, что задача~\eqref{eq:problem_statement} обладает $\epsilon$-острым минимумом $(\epsilon>0)$, если для фиксированного $\alpha>0$ верна альтернатива
    \begin{itemize}
        \item при $g(x)>\epsilon$
        $$
        g(x) \geq \alpha \min _{x_{*} \in X_{*}}\left\|x-x_{*}\right\|_{2} \quad \forall x \in Q,
        $$
        \item при $g(x) \leq \epsilon$
        $$
        f(x)-f^{*} \geq \alpha \min _{x_{*} \in X_{*}}\left\|x-x_{*}\right\|_{2} \quad \text { или } \quad f(x)-f^{*} \leq \epsilon \quad \forall x \in Q .
        $$
    \end{itemize}
\end{fed}

\begin{exl}
    Рассмотрим задачу минимизации сильно выпуклой функции
    \begin{gather}\label{ex_dist}
    \min _{g(x) \leq 0, x \in Q} f(x),
    \end{gather}
    где $Q$~--- выпуклое ограниченное замкнутое множество, а $g(x) = \text{dist} (x, Q) = \min\limits_{y\in Q}\|x-y\|_2$.\\
    Действительно, задача~\eqref{ex_dist} обладает $\epsilon$-острым минимумом: т.~к. $Q$~--- ограниченное множество, то diam$(Q) \leq C$ для некоторого $C>0$, поэтому в случае $g(x) = \text{dist}(x, Q) > \epsilon$ воспользуемся неравенством треугольника
    $$
    \min_{x_*\in X_*}\|x-x_*\|_2 \leq \text{dist}(x, Q) + \text{diam}(Q) \leq g(x) + C,
    $$
    поэтому
    $$
    g(x) \geq \min_{x_*\in X_*}\|x-x_*\|_2 - C = \alpha \min_{x_*\in X_*}\|x-x_*\|_2.
    $$
    При $g(x) \leq \epsilon$ всё так же потребуем от целевой функции условие острого минимума:
    $$
    f(x)-f^{*} \geq \alpha \min _{x_{*} \in X_{*}}\left\|x-x_{*}\right\|_{2} \quad \text { или } \quad f(x)-f^{*} \leq \epsilon,
    $$
    для некоторого $\alpha>0.$
\end{exl}


Покажем, что c использованием некоторого варианта шага Б.\,Т.~Поляка можно достичь сходимости следующего метода со скоростью геометрической прогрессии к $\varepsilon$-точному решению  поставленной задачи.

\begin{algorithm}[!ht]
\caption{Схема с переключениями по продуктивным и непродуктивным шагам для липшицевых задач с ограничениями-неравенствами c $\epsilon$-острым минимумом.}
\label{alg:eps_sharp_min}
\begin{algorithmic}[1]
    \REQUIRE $x_0: \, \min\limits_{x_* \in X_*} \norm*{x_0 - x_*}_2 \leq \frac{\alpha \gamma}{\mu}$, $\gamma \in(0; 1)$, $\epsilon > 0, M_f > 0$.
    \STATE $I \gets \emptyset$, $J \gets \emptyset$, $k \gets 0$;
    \REPEAT
        \IF{$g(x_k) \leq \epsilon$}
            \STATE $\displaystyle h_{k}^{f} = \frac{f\left(x_{k}\right)-\overline{f}}{M_f\left\|\nabla f\left(x_{k}\right)\right\|_{2}}$;
            \STATE $\displaystyle x_{k+1} = \operatorname{Pr}_{Q}\left\{x_{k}-h_{k}^{f} \nabla f\left(x_{k}\right)\right\}$; \hspace{1cm} <<продуктивные шаги>>
            \STATE $k \to I$;
        \ELSE
            \STATE $\displaystyle h_{k}^{g} = \frac{g\left(x_{k}\right)}{\left\|\nabla g\left(x_{k}\right)\right\|^2_{2}}$;
            \STATE $\displaystyle x_{k+1} = \operatorname{Pr}_{Q}\left\{x_{k}-h_{k}^{g} \nabla g\left(x_{k}\right)\right\}$; \hspace{1cm} <<непродуктивные шаги>>
            \STATE $k \to J$;
        \ENDIF
        \STATE $k \gets k + 1$;
    \UNTIL $k = N$.
\end{algorithmic}
\end{algorithm}

\begin{teo}\label{teo:eps_sharp_min}
    Пусть функция $f$ квазивыпукла и $M_f$-липшицева ($M_f>0$), a $g$ $\mu$-слабо выпукла и задача имеет $\epsilon$-острый минимум ($\alpha>0$, $\mu>0$). Пусть $\overline{f}$ такое, что $f(x) - \overline{f} = c(x)\left(f(x)-f^*\right)$,
    где $c(x) \in [C, 2-C]$ для некоторого $C \in (0, 1)$. Начальная точка $x_{0}$ такова, что $\underset{x_{*} \in X_{*}}{\min}\left\|x_{0}-x_{*}\right\|_{2} \leq \frac{\alpha \gamma}{\mu}$ для некоторого фиксированного $\gamma \in(0 ; 1)$. Тогда после $k+1$ итераций алгоритма~\ref{alg:eps_sharp_min} справедлива следующая альтернатива: либо достигнуто $\epsilon$-точное решение задачи, т.~е.
    $$
    f\left(x_{k+1}\right)-f^{*} \leq \epsilon, \quad g\left(x_{k+1}\right) \leq \epsilon,
    $$
    или верно следующее неравенство:
    $$
    \min _{x_{*} \in X_{*}}\left\|x_{k+1}-x_{*}\right\|_{2}^{2} \leq \left(1-\frac{\alpha^{2}}{M_f^2}\left( 2C-C^2\right)\right)^{\left|I\right|} \prod_{i \in J}\left(1-\frac{\left(1-\gamma_{i}\right) \alpha^{2}}{\left\|\nabla g\left(x_{i}\right)\right\|_{2}^{2}}\right) \min _{x_{*} \in X_{*}}\left\|x_{0}-x_{*}\right\|_{2}^{2},
    $$
    где последовательность $\left\{\gamma_{i}\right\}_{i=0}^{k}$ такова, что $\gamma_{0}:=\gamma \in(0 ; 1)$, а далее, при $i=0,1,2, \ldots, k$, определим
    $$
    \gamma_{i+1}:=\gamma_{i} \sqrt{1-\frac{\alpha^{2}(1-\gamma_{i})}{\left\|\nabla g\left(x_{i}\right)\right\|_{2}^{2}}}
    $$
    для всякого непродуктивного шага $(i \in J)$ и 
    $$
    \gamma_{i+1} := \gamma_{i}\sqrt{1 - \frac{\alpha^2}{M_f^2}\left(2C-C^2\right)}
    $$
    для всякого продуктивного шага $\left(i\in I\right)$.
\end{teo}
\proof*
    1) Для продуктивных шагов по лемме~\ref{lem:proj} справедливо следующее неравенство:
    $$
    2 h_{k}^{f}\left\langle\nabla f\left(x_{k}\right), x_{k}-x_{*}\right\rangle \leq \left(h_{k}^{f}\right)^{2}\left\|\nabla f\left(x_{k}\right)\right\|_{2}^{2}+\left\|x_{k}-x_{*}\right\|_{2}^{2}-\left\|x_{k+1}-x_{*}\right\|_{2}^{2} \quad \forall k \in I.
    $$
    Поэтому имеем
    $$
    \begin{aligned}
        \min _{x_{*} \in X_{*}}\left\|x_{k+1}-x_{*}\right\|_{2}^{2} & \leq 
         \min _{x_{*} \in X_{*}}\left\|x_{k}-x_{*}\right\|_{2}^{2}+\left(h_{k}^{f}\right)^{2}\left\|\nabla f\left(x_{k}\right)\right\|_{2}^{2}-2 h_{k}^{f}\left\langle\nabla f\left(x_{k}\right), x_{k}-x_{*}\right\rangle 
        \\& 
        = \min _{x_{*} \in X_{*}}\left\|x_{k}-x_{*}\right\|_{2}^{2}+\frac{\left(f\left(x_{k}\right)-\overline{f}\right)^{2}}{M_f^{2}}-\frac{2\left(f\left(x_{k}\right)-\overline{f}\right)}{M_f}\left\langle\frac{\nabla f\left(x_{k}\right)}{\left\|\nabla f\left(x_{k}\right)\right\|_{2}}, x_{k}-x_{*}\right\rangle 
        \\&
        = \min _{x_{*} \in X_{*}}\left\|x_{k}-x_{*}\right\|_{2}^{2}+\left(\frac{f\left(x_{k}\right)-\overline{f}}{M_f}\right)^{2}-\frac{2\left(f\left(x_{k}\right)-\overline{f}\right)}{M_f} v_{f}\left(x_{k}, x_{*}\right) 
        \\&
        \leq \min _{x_{*} \in X_{*}}\left\|x_{k}-x_{*}\right\|_{2}^{2} + c^2(x_k)\frac{(f(x_k)-f^*)^2}{M_f^2} - 2 \cdot\frac{f(x_k) - f^*}{M_f^2}\cdot c(x_k)\left(f(x_k) - f^*\right) 
        \\&
        = \min _{x_{*} \in X_{*}}\left\|x_{k}-x_{*}\right\|_{2}^{2} - \left(\frac{f(x_k)-f^*}{M_f}\right)^2\left(2c(x_k)-c^2(x_k)\right) 
        \\&
        \leq \min _{x_{*} \in X_{*}}\left\|x_{k}-x_{*}\right\|_{2}^{2} \left(1 - \frac{\alpha^2}{M_f^2}\left(2c(x_k)-c^2(x_k)\right)\right) 
        \\&
        \leq\min _{x_{*} \in X_{*}}\left\|x_{k}-x_{*}\right\|_{2}^{2} \left(1 - \frac{\alpha^2}{M_f^2}\left(2C-C^2\right)\right),
    \end{aligned}
    $$
    так как $v_{f}\left(x_{k}, x_{*}\right) \geq \frac{f\left(x_{k}\right)-f^{*}}{M_f}$. \\\\
    2) Для непродуктивных шагов при $g\left(x_{k}\right)>\epsilon$ верно $g\left(x_{k}\right) \geq \alpha \underset{x_{*} \in X_{*}}{\min} \left\|x_{k}-x_{*}\right\|_{2}$. Тогда по лемме~\ref{lem:proj} для ближайшего к $x_k$ точного решения $x_{*} \in X_{*}$ ввиду $g\left(x_{*}\right) \leq 0$, имеем:
    $$
    \begin{aligned}
    \left\|x_{k+1}-x_{*}\right\|_{2}^{2} & \leq
    \left\|x_{k}-x_{*}\right\|_{2}^{2}+\frac{2 g\left(x_{k}\right)\left\langle\nabla g\left(x_{k}\right), x_{*}-x_{k}\right\rangle}{\left\|\nabla g\left(x_{k}\right)\right\|_{2}^{2}}+\left(\frac{g\left(x_{k}\right)}{\left\|\nabla g\left(x_{k}\right)\right\|_{2}^{2}}\right)^{2}\left\|\nabla g\left(x_{k}\right)\right\|_{2}^{2} 
    \\&
    \leq\left\|x_{k}-x_{*}\right\|_{2}^{2}+\frac{2 g\left(x_{k}\right)}{\left\|\nabla g\left(x_{k}\right)\right\|_{2}^{2}}\left(g\left(x_{*}\right)-g\left(x_{k}\right)+\frac{\mu}{2}\left\|x_{k}-x_{*}\right\|_{2}^{2}\right)+\frac{\left(g\left(x_{k}\right)\right)^{2}}{\left\|\nabla g\left(x_{k}\right)\right\|_{2}^{2}} 
    \\&
    \leq\left\|x_{k}-x_{*}\right\|_{2}^{2}+\frac{g\left(x_{k}\right) \mu}{\left\|\nabla g\left(x_{k}\right)\right\|_{2}^{2}}\left\|x_{k}-x_{*}\right\|_{2}^{2}-\frac{2\left(g\left(x_{k}\right)\right)^{2}}{\left\|\nabla g\left(x_{k}\right)\right\|_{2}^{2}}+\frac{\left(g\left(x_{k}\right)\right)^{2}}{\left\|\nabla g\left(x_{k}\right)\right\|_{2}^{2}}
    \\&
    =\left\|x_{k}-x_{*}\right\|_{2}^{2}+\frac{g\left(x_{k}\right)}{\left\|\nabla g\left(x_{k}\right)\right\|_{2}^{2}}\left(\mu\left\|x_{k}-x_{*}\right\|_{2}^{2}-g\left(x_{k}\right)\right) 
    \\&
    \leq\left\|x_{k}-x_{*}\right\|_{2}^{2}+\frac{g\left(x_{k}\right)}{\left\|\nabla g\left(x_{k}\right)\right\|_{2}^{2}}\left(\mu\left\|x_{k}-x_{*}\right\|_{2}^{2}-\alpha\left\|x_{k}-x_{*}\right\|_{2}\right) 
    \\&
    \leq\left\|x_{k}-x_{*}\right\|_{2}^{2}+\frac{g\left(x_{k}\right)}{\left\|\nabla g\left(x_{k}\right)\right\|_{2}^{2}}\left(\alpha \gamma_{k}-\alpha\right)\left\|x_{k}-x_{*}\right\|_{2} \leq\left(1-\frac{\left(1-\gamma_{k}\right) \alpha^{2}}{\left\|\nabla g\left(x_{k}\right)\right\|_{2}^{2}}\right)\left\|x_{k}-x_{*}\right\|_{2}^{2} .
    \end{aligned}
    $$
    Тогда ясно, что
    $$
    \min _{x_{*} \in X_{*}}\left\|x_{k+1}-x_{*}\right\|_{2}^{2} \leq \left(1-\frac{\alpha^{2}}{M_f^2}\left(2C-C^2\right)\right)^{\left|I\right|} \prod_{i \in J}\left(1-\frac{\left(1-\gamma_{i}\right) \alpha^{2}}{\left\|\nabla g\left(x_{i}\right)\right\|_{2}^{2}}\right) \min _{x_{*} \in X_{*}}\left\|x_{0}-x_{*}\right\|_{2}^{2},
    $$
    что и требовалось. При этом на <<непродуктивном шаге>> не может быть ситуации $\nabla g(x_k) = 0$ ввиду предложенного варианта выбора начальной точки $x_0$ (\cite{Stonyakin2023}, замечание 4).
\qed

Если рассматривать задачи с выпуклыми ограничениями ($\mu = 0$), то пропадает необходимость ограничивать выбор начальной точки $x_0$ в алгоритме~\ref{alg:eps_sharp_min}. Сформулируем этот результат.

\begin{cor}
    Пусть $f$ квазивыпукла и $M_f$-липшицева ($M_f>0$), a $g$ выпукла и задача имеет $\epsilon$-острый минимум ($\alpha>0$). Тогда после $k+1$ итераций алгоритма~\ref{alg:eps_sharp_min} справедлива следующая альтернатива: либо достигнуто $\epsilon$-точное решение задачи, т.~е.
    $$
    f\left(x_{k+1}\right)-f^{*} \leq \epsilon, \quad g\left(x_{k+1}\right) \leq \epsilon,
    $$
    или верно следующее неравенство:
    $$
    \min _{x_{*} \in X_{*}}\left\|x_{k+1}-x_{*}\right\|_{2}^{2} \leq \left(1-\frac{\alpha^{2}}{M_f^2}\left(2C-C^2\right)\right)^{\left|I\right|} \prod_{i\in J}\left(1-\frac{\alpha^{2}}{\left\|\nabla g\left(x_{i}\right)\right\|_{2}^{2}}\right) \min _{x_{*} \in X_{*}}\left\|x_{0}-x_{*}\right\|_{2}^{2}.
    $$
\end{cor}
\proof*
    Аналогично теореме~\ref{teo:eps_sharp_min} с $\mu = 0$.
\qed

\begin{rem}
    Отметим, что равномерная ограниченность норм $\left\|\nabla g\left(x_{k}\right)\right\|_{2} \leq M_g,$ $ \forall k \in \overline{0, N}$, где $M_g > 0$ ,  влечет в выпуклом случае оценку
    $$
    \begin{aligned}
    \min _{x_{*} \in X_{*}}\left\|x_{N}-x_{*}\right\|_{2}^{2} & \leq
    \left(1-\frac{\alpha^{2}}{M_f^2}\left(2C-C^2\right)\right)^{\left|I\right|} \left(1-\frac{\alpha^{2}}{M_{g}^{2}}\right)^{\left|J\right|} \min _{x_{*} \in X_{*}}\left\|x_{0}-x_{*}\right\|_{2}^{2}
    \\& 
    \leq \left(1-\frac{\alpha^2}{M^2}\left(2C-C^2\right)\right)^{N}\min _{x_{*} \in X_{*}}\left\|x_{0}-x_{*}\right\|_{2}^{2},
    \end{aligned}
    $$
    где $M = \max\{M_f, M_g\}$. Это, в свою очередь, указывает на сходимость со скоростью геометрической прогрессии, причем при $C=1$ получаем $\overline{f} = f^*$ и наибольшую скорость сходимости.
\end{rem}

\begin{rem}\label{rem:G}
Также шаг Б.\,Т.~Поляка позволяет доказать, что невязка по аргументу не возрастает по мере увеличения количества итераций. А значит, можно считать, что для любого $k = 0, 1, 2, \ldots$ выполняется $x_k \in G$, где $G$ — некоторое ограниченное множество. Поэтому достаточно требовать ограниченность норм градиента и, как следствие, липшицевость функции не на всем пространстве, а лишь на множестве $G$.
\end{rem}

\begin{rem}
    В некоторых задачах возникает условие \textit{слабого острого минимума} 
    \begin{gather}
        f(x) - f^* \geq C \min_{x_*\in X_*}\|x-x_*\|_2^p \quad \forall x \in Q,
    \end{gather} для некоторого фиксированного $p\in[1, +\infty)$ и некоторой постоянной $C>0$. Тогда при достаточно большом расстоянии от текущей точки до множества решений $X_*$ выполняется $$f(x) - f^* \geq C \min_{x_*\in X_*}\|x-x_*\|_2^p \geq C_1 \min_{x_*\in X_*}\|x-x_*\|_2,$$ где $C_1$~--- некоторая постоянная. Таким образом, описанная методика позволит достигнуть определенного уровня близости к точному решению по аргументу.
\end{rem}

\begin{rem}    
Пусть теперь задача ставится следующим образом:
\begin{gather}\label{many_constraints}
    \min _{g(x) \leq 0, x \in Q} f(x),
\end{gather}
где $Q$~--- выпуклое замкнутое подмножество $\mathbb{R}^n$, а $g(x) = \max\limits_{i \in \overline{1, m}} g_i(x)$. Тогда следует переформулировать определение $\epsilon$-острого минимума задачи~\eqref{many_constraints}.
\begin{fed}
    Будем говорить, что задача~\eqref{many_constraints} обладает $\epsilon$-острым минимумом $(\epsilon>0)$, если существует набор $\alpha_{g_i}>0$ для $i \in \overline{1, m}$ и $\alpha_f > 0$, для которых верна альтернатива
    \begin{itemize}
        \item при $g_i(x)>\epsilon$,
        $$
        g_i(x) \geq \alpha_{g_i} \min _{x_{*} \in X_{*}}\left\|x-x_{*}\right\|_{2} \quad \forall x \in Q,
        $$
        \item при $g_i(x) \leq \epsilon$,
        $$
        f(x)-f^{*} \geq \alpha_f \min _{x_{*} \in X_{*}}\left\|x-x_{*}\right\|_{2} \quad \text { или } \quad f(x)-f^{*} \leq \epsilon \quad \forall x \in Q .
        $$
    \end{itemize}
\end{fed}
Тогда множество продуктивных шагов алгоритма обозначим $I = \{k\in \mathbb{N} \mid \forall i \in \overline{1,m} \hookrightarrow g(x_k)\leq \epsilon\}$, непродуктивных~--- $J = \{k\in \mathbb{N} \mid \exists i: g_i(x_k)>\epsilon\}$. То есть пропадает необходимость считать все субградиенты $\nabla g_i(x)$, вместо этого работая с первым нарушенным ограничением. При этом на непродуктивных шагах в качестве параметров ограничений ($\alpha, M, \mu$ и т.~д.) можно выбрать наиболее удачные с точки зрения сходимости.
\end{rem}

\paragraph{4. Cубградиентный метод с шагом типа Б.\,Т.~Поляка для липшицевых задач с ограничениями-неравенствами с вариантом острого минимума, не зависящим от ожидаемой точности решения задач}

В этом разделе мы рассмотрим другой способ задания острого минимума, предложенный в~\cite{Lin2020}. 
\begin{fed}
    Будем говорить, что задача~\eqref{eq:problem_statement} обладает <<условным>> острым минимумом, если существует такое $\alpha > 0$, для которого выполнено
    \begin{gather}\label{eq:sharp_min_without_eps}
        \max\cbr*{f(x) - f^*, g(x)} \geq \alpha \min_{x_* \in X_*} \norm*{x - x_*}_2.
    \end{gather}
\end{fed}
В недавней работе~\cite{Ablaev2023} для липшицевых задач с <<условным>> острым минимумом доказана сходимость со скоростью геометрической прогрессии рестартованных по параметру острого минимума субградиентных методов. Недостаток этого подхода~--- необходимость знать оценку такого параметра, что довольно затруднительно для многих реально возникающих задач на практике. Мы поставим задачу получить результат о линейной скорости сходимости (хотя бы в окрестность множества решений) для варианта субградиентного метода, не предполагающего для реализации знания величины этого параметра.
В отличие от рассмотренного в предыдущем разделе $\epsilon$-острого минимума, данный способ задания острого минимума не зависит от параметра $\epsilon$, который, как видно из теоремы~\ref{teo:eps_sharp_min}, ограничивает точность, которую может гарантировать алгоритм~\ref{alg:eps_sharp_min}. Предлагаемый же алгоритм~\ref{alg:sharp_min_without_eps} для задачи с острым минимумом~\eqref{eq:sharp_min_without_eps}, как показано далее в этом разделе, в теореме~\ref{teo:sharp_min_without_eps}, имеет линейную скорость сходимости и при достаточно большом числе итераций может позволить достичь любого заданного качества решения по аргументу в зависимости от точности информации о минимальном значении $f^*$.

\begin{algorithm}[!ht]
\caption{Схема с переключениями по продуктивным и непродуктивным шагам для липшицевых задач с ограничениями-неравенствами с вариантом острого минимума, не зависящим от ожидаемой точности решения задач. }
\label{alg:sharp_min_without_eps}
\begin{algorithmic}[1]
    \REQUIRE $x_0: \, \min\limits_{x_* \in X_*} \norm*{x_0 - x_*}_2 \leq \frac{C \alpha \gamma_0}{\mu}$, $\gamma_0 \in (0; 1), M_f > 0, N > 0$.
    \STATE $I \gets \emptyset$, $J \gets \emptyset$, $k \gets 0$;
    \REPEAT
        \IF{$f(x_k) - \overline{f} \geq g(x_k)$}
            \STATE $\displaystyle h_k^f = \frac{f(x_k) -\overline{f}}{M_f \norm*{\nabla f(x_k)}_2}$;
            \STATE $\displaystyle x_{k+1} = \operatorname{Pr}_{Q} \cbr*{x_k - h_k^f \nabla f(x_k)}$; \hspace{1cm} <<продуктивные шаги>>
            \STATE $k \to I$;
        \ELSE
            \STATE $\displaystyle h_k^g = \frac{g(x_k)}{\norm*{\nabla g(x_k)}_2^2}$;
            \STATE $\displaystyle x_{k+1} = \operatorname{Pr}_{Q} \cbr*{x_k - h_k^g \nabla g(x_k)}$; \hspace{1cm} <<непродуктивные шаги>>
            \STATE $k \to J$;
        \ENDIF
        \STATE $k \gets k + 1$;
    \UNTIL $k = N$.
\end{algorithmic}
\end{algorithm}

Проведем анализ алгоритма~\ref{alg:sharp_min_without_eps} для задачи~\eqref{eq:problem_statement} с учётом всех предположений (отсутствие стационарных точек в $Q \setminus X_*$, наличие острого минимума~\eqref{eq:sharp_min_without_eps} и доступ к $\overline{f}$).

Начнём с анализа <<продуктивных шагов>>. Возможны два случая:
\begin{itemize}
    \item $f(x_k) - f^* \geq g(x_k)$: в силу условия острого минимума~\eqref{eq:sharp_min_without_eps} и условия~\eqref{eq:f_plus} на $\overline{f}$ выполнено
        \begin{gather}
            f(x_k) - \overline{f} = c(x_k) \rbr*{f(x_k) - f^*} \geq c(x_k) \alpha \min_{x_* \in X_*} \norm*{x_k - x_*}_2;
        \end{gather}
    \item $f(x_k) - f^* \leq g(x_k)$: тогда, так как шаг <<продуктивный>> и в силу условия острого минимума~\eqref{eq:sharp_min_without_eps}, выполнено
        \begin{gather}
            f(x_k) - \overline{f} \geq g(x_k) \geq \alpha \min_{x_* \in X_*} \norm*{x_k - x_*}_2.
        \end{gather}
\end{itemize}
То есть
\begin{gather}\label{eq:prod_init}
    f(x_k) - \overline{f} \geq \min\cbr*{1, c(x_k)} \alpha \min_{x_* \in X_*} \norm*{x_k - x_*}_2.
\end{gather}

Рассмотрим ближайшее $x_* \in X_*$ к $x^k$ и, применив лемму~\ref{lem:proj} к $x_{k+1}$ на <<продуктивном шаге>>, получим
\begin{gather}
    \begin{aligned}
        \norm*{x_{k+1} - x_*}_2^2
        &\leq \norm*{x_k - x_*}_2^2 - 2 h_k^f \abr*{\nabla f(x_k), x_k - x_*} + \rbr*{h_k^f}^2 \norm*{\nabla f(x_k)}_2^2 \\
        &= \norm*{x_k - x_*}_2^2 - 2 \frac{f(x_k) -\overline{f}}{M_f} v_f(x_k, x_*) + \rbr*{\frac{f(x_k) -\overline{f}}{M_f}}^2.
    \end{aligned}
\end{gather}
Воспользовавшись леммой~\ref{lem:v_f} и условием~\eqref{eq:f_plus} на $\overline{f}$, получаем
\begin{gather}
    \norm*{x_{k+1} - x_*}_2^2 \leq \norm*{x_k - x_*}_2^2 - \frac{2 - c(x_k)}{c(x_k)} \rbr*{\frac{f(x_k) -\overline{f}}{M_f}}^2.
\end{gather}
Наконец, применив~\eqref{eq:prod_init},
\begin{gather}
\norm*{x_{k+1} - x_*}_2^2
\leq \rbr*{1 - \frac{2 - c(x_k)}{c(x_k)} \rbr*{\min\cbr*{1, c(x_k)}}^2 \frac{\alpha^2}{M_f^2}} \norm*{x_k - x_*}_2^2 
\leq \rbr*{1 - \frac{C}{2 - C} \frac{\alpha^2}{M_f^2}} \norm*{x_k - x_*}_2^2.
\end{gather}
Следовательно,
\begin{gather}\label{eq:prod}
    \min_{x_* \in X_*} \norm*{x_{k+1} - x_*}_2^2 \leq \rbr*{1 - \frac{C}{2 - C} \frac{\alpha^2}{M_f^2}} \min_{x_* \in X_*} \norm*{x_k - x_*}_2^2.
\end{gather}

Теперь перейдём к анализу <<непродуктивного шага>>. Действовать будем по аналогичной схеме. Возможны два случая:
\begin{itemize}
    \item $f(x_k) - f^* \geq g(x_k)$: в силу условия острого минимума~\eqref{eq:sharp_min_without_eps}, условия~\eqref{eq:f_plus} на $\overline{f}$ и того, что шаг <<непродуктивный>>, выполнено
        \begin{gather}
            g(x_k) \geq f(x_k) - \overline{f} = c(x_k) \rbr*{f(x_k) - f^*} \geq c(x_k) \alpha \min_{x_* \in X_*} \norm*{x_k - x_*}_2 \geq C \alpha \min_{x_* \in X_*} \norm*{x_k - x_*}_2;
        \end{gather}
    \item $f(x_k) - f^* \leq g(x_k)$: тогда, в силу условия острого минимума~\eqref{eq:sharp_min_without_eps}, выполнено
        \begin{gather}
            g(x_k) \geq \alpha \min_{x_* \in X_*} \norm*{x_k - x_*}_2.
        \end{gather}
\end{itemize}
То есть, так как $C \leq 1$,
\begin{gather}\label{eq:notprod_init}
    g(x_k) \geq C \alpha \min_{x_* \in X_*} \norm*{x_k - x_*}_2.
\end{gather}

Снова взяв ближайшее $x_* \in X_*$ к $x^k$ и воспользовавшись леммой~\ref{lem:proj} для $x_{k+1}$ с <<непродуктивного шага>>, получим
\begin{gather}
    \begin{aligned}
        \norm*{x_{k+1} - x_*}_2^2
        &\leq \norm*{x_k - x_*}_2^2 - 2 h_k^g \abr*{\nabla g(x_k), x_k - x_*} + \rbr*{h_k^g}^2 \norm*{\nabla g(x_k)}_2^2 \\
        &= \norm*{x_k - x_*}_2^2 - 2 \frac{g(x_k)}{\norm*{\nabla g(x_k)}_2^2} \abr*{\nabla g(x_k), x_k - x_*} + \frac{\rbr*{g(x_k)}^2}{\norm*{\nabla g(x_k)}_2^2}.
    \end{aligned}
\end{gather}
В силу $\mu$-слабой выпуклости $g$ имеем
\begin{gather}
    \begin{aligned}
        \norm*{x_{k+1} - x_*}_2^2
        &\leq \norm*{x_k - x_*}_2^2 + 2 \frac{g(x_k)}{\norm*{\nabla g(x_k)}_2^2} \rbr*{g(x_*) - g(x_k) + \frac{\mu}{2} \norm*{x_k - x_*}_2^2} + \frac{\rbr*{g(x_k)}^2}{\norm*{\nabla g(x_k)}_2^2} \\
        &\leq \norm*{x_k - x_*}_2^2 - \frac{g(x_k)}{\norm*{\nabla g(x_k)}_2^2} \rbr*{g(x_k) - \mu \norm*{x_k - x_*}_2^2},
    \end{aligned}
\end{gather}
где в последнем неравенстве учтено, что $g(x_*) \leq 0$. Далее, дважды воспользовавшись неравенством~\eqref{eq:notprod_init}, получим
\begin{gather}
    \begin{aligned}
        \norm*{x_{k+1} - x_*}_2^2
        &\leq \norm*{x_k - x_*}_2^2 - \frac{g(x_k)}{\norm*{\nabla g(x_k)}_2^2} \rbr*{C \alpha \norm*{x_k - x_*}_2 - \mu \norm*{x_k - x_*}_2^2} \\
        &\leq \norm*{x_k - x_*}_2^2 - \frac{g(x_k)}{\norm*{\nabla g(x_k)}_2^2} (1 - \gamma_k) C \alpha \norm*{x_k - x_*}_2 \\
        &\leq \rbr*{1 - \frac{(1 - \gamma_k) C^2 \alpha^2}{\norm*{\nabla g(x_k)}_2^2}} \norm*{x_k - x_*}_2^2,
    \end{aligned}
\end{gather}
где последовательность $\cbr*{\gamma_{i+1}}_{i=0}^{k-1}$ задается следующим образом:
\begin{gather}\label{eq:gamma_k_without_eps}
    \gamma_{i+1} =  \begin{cases}
        \gamma_i \sqrt{1 - \frac{C}{2 - C} \frac{\alpha^2}{M_f^2}}, \quad &i \in I;\\
        \gamma_i \sqrt{1 - \frac{(1 - \gamma_i) C^2 \alpha^2}{\norm*{\nabla g(x_i)}_2^2}}, \quad &i \in J,
    \end{cases}
\end{gather}
гарантирующим, что $\min\limits_{x_* \in X_*} \norm*{x_k - x_*}_2 \leq \frac{C \alpha \gamma_k}{\mu}$. Следовательно,
\begin{gather}\label{eq:notprod}
    \min_{x_* \in X_*} \norm*{x_{k+1} - x_*}_2^2 \leq \rbr*{1 - C^2 \frac{(1 - \gamma_k) \alpha^2}{\norm*{\nabla g(x_k)}_2^2}} \min_{x_* \in X_*} \norm*{x_k - x_*}_2^2.
\end{gather}

Осталось показать, что на <<непродуктивном шаге>> не может быть ситуации $\nabla g(x_k) = 0$. Предположим противное. Тогда, в силу $\mu$-слабой выпуклости функции $g$, для любого $x_* \in X_*$
\begin{gather}\label{eq:tmp}
    g(x_*) \geq g(x_k) + \abr*{\nabla g(x_k), x_* - x_k} - \frac{\mu}{2} \norm*{x_* - x_k}_2^2 = g(x_k) - \frac{\mu}{2} \norm*{x_* - x_k}_2^2.
\end{gather}
Так как любое $x_* \in X_*$ является решением, то $g(x_*) \leq 0$. Следовательно, заметив, что неравенство~\eqref{eq:tmp} верно для любого $x_* \in X_*$, и воспользовавшись неравенством~\eqref{eq:notprod_init}, имеем
\begin{gather}
    C \alpha \min_{x_* \in X_*} \norm*{x_k - x_*}_2 \leq g(x_k) \leq \frac{\mu}{2} \min_{x_* \in X_*} \norm*{x_* - x_k}_2^2.
\end{gather}
То есть $\min\limits_{x_* \in X_*} \norm*{x_* - x_k}_2 \geq \frac{2 C \alpha}{\mu}$, что неверно, так как $\min\limits_{x_* \in X_*} \norm*{x_0 - x_*}_2 \leq \frac{C \alpha \gamma_0}{\mu}$, где $\gamma_0 \in (0; 1)$, и, в силу оценок~\eqref{eq:prod} и \eqref{eq:notprod}, расстояния от каждой следующей получаемой алгоритмом точки до множества решений меньше, чем от каждой предыдущей. 

Таким образом, верна следующая
\begin{teo}\label{teo:sharp_min_without_eps}
Пусть $Q \subset \RR^n$~--- замкнутое выпуклое множество, $f: Q \longrightarrow \RR$~--- $M_f$-липшицевая квазивыпуклая функция, $g: Q \longrightarrow \RR$~--- $\mu$-слабо выпуклая функция, а начальная точка $x_0$ такова, что $\min\limits_{x_* \in X_*} \norm*{x_0 - x_*}_2 \leq \frac{C \alpha \gamma_0}{\mu}$, где $\gamma_0 \in (0; 1)$. Тогда верно следующее неравенство для выходной точки $x_N$ алгоритма~\ref{alg:sharp_min_without_eps} после $N$ итераций:
\begin{gather}
    \min_{x_* \in X_*} \norm*{x_N - x_*}_2^2 \leq \rbr*{1 - \frac{C}{2 - C} \frac{\alpha^2}{M_f^2}}^{\cardi*{I}} \prod_{k \in J} \rbr*{1 - C^2 \frac{(1 - \gamma_k) \alpha^2}{\norm*{\nabla g(x_k)}_2^2}} \min_{x_* \in X_*} \norm*{x_0 - x_*}_2^2,
\end{gather}
где последовательность $\cbr*{\gamma_{i+1}}_{i=0}^{N-1}$ задается согласно~\eqref{eq:gamma_k_without_eps}.
\end{teo}

\begin{cor}
    Если в условиях теоремы~\ref{teo:sharp_min_without_eps} дополнительно предположить выпуклость функции $g$, или, иными словами, что $\mu = 0$, то будет верно следующее неравенство для выходной точки $x_N$ алгоритма~\ref{alg:sharp_min_without_eps} после $N$ итераций:
    \begin{gather}
        \min_{x_* \in X_*} \norm*{x_N - x_*}_2^2 \leq \rbr*{1 - \frac{C}{2 - C} \frac{\alpha^2}{M_f^2}}^{\cardi*{I}} \prod_{k \in J} \rbr*{1 - C^2 \frac{\alpha^2}{\norm*{\nabla g(x_k)}_2^2}} \min_{x_* \in X_*} \norm*{x_0 - x_*}_2^2.
    \end{gather}
    При этом выбор начальной точки ничем не ограничен.
\end{cor}
\proof*
    Аналогично приведенным выше рассуждениям при $\mu = 0$.
\qed

\begin{cor}\label{cor:21}
    Если, помимо выпуклости функции $g$, в условиях теоремы~\ref{teo:sharp_min_without_eps} предположить ограниченность норм $\|\nabla g(x_k)\|_2 \leq M_g$ для всех $k = \overline{0,N-1}$ (например, в случае $M_g$-липшицевости функции $g$), то будет верно следующее неравенство для выходной точки $x_N$ алгоритма~\ref{alg:sharp_min_without_eps} после $N$ итераций:
    \begin{gather}
        \min_{x_* \in X_*} \norm*{x_N - x_*}_2^2 \leq \rbr*{1 - C^2 \frac{\alpha^2}{M^2}}^{N} \min_{x_* \in X_*} \norm*{x_0 - x_*}_2^2,
    \end{gather}
    где $M = \max\cbr*{M_f, M_g}$.
\end{cor}
\proof*
    Аналогично приведенным выше рассуждениям с учетом того, что $\frac{C}{2 - C} \geq C^2$, так как $C \in (0; 1]$.
\qed

\begin{rem}
    Аналогично замечанию~\ref{rem:G} для алгоритма~\ref{alg:eps_sharp_min} липшицевость функции $f$, а также ограниченность норм в следствии~\ref{cor:21} можно требовать не на всём множестве Q, а лишь на некотором ограниченном подмножестве $G$.
\end{rem}

\paragraph{5. Вычислительные эксперименты}

В настоящем разделе для иллюстрации работоспособности предложенных выше алгоритмов~\ref{alg:eps_sharp_min} и~\ref{alg:sharp_min_without_eps} приведём некоторые результаты вычислительных экспериментов для четырёх примеров в сравнении с алгоритмом 3 из \cite{Ablaev2023}, который является одним из недавно предложенных алгоритмов для класса квазивыпуклых задач оптимизации с функциональными ограничениями-неравенствами.

\begin{exl}[задача геометрического программирования \cite{boyd2007tutorial}]\label{ex_geometric_prob}  
Рассмотрим целевую функцию вида
    \begin{gather}\label{obj_lp}
        f(x) = \sqrt[p]{|x_1|^p + |x_2|^p + \ldots + |x_n|^p }, \quad p \geq 1; 
    \end{gather}
и следующую функцию ограничений:
    \begin{gather}\label{cons_gemetric_prob}
        g(x) = \max_{i \in \overline{1, m}} \left\{ g_i(x) = a_i x_1^{\alpha_{i1} } x_2^{\alpha_{i2} } \cdots x_n^{\alpha_{in} } - b_i, \; x_j >0 \, \forall j \in \overline{1, n} \right\}, 
    \end{gather}
где $a_i > 0, b_i \in \mathbb{R},  (\alpha_{i1}, \alpha_{i2}, \ldots, \alpha_{in}) \in \mathbb{R}^n$ и каждая функция $g_i(x)\,  (i \in \overline{1, m})$ является позиномиальной функцией.
\end{exl}

\begin{exl}\label{ex_quasi}
Предположим, что $a, b \in \mathbb{R}^n$, и определим следующую целевую функцию 
    \begin{gather}\label{obj_quasi}
        f(x) = \frac{\|x-a\|_2}{\|x-b\|_2}, 
    \end{gather}
то есть отношение евклидова расстояния до $a$ к расстоянию до $b$. Функция $f$ квазивыпукла на полупространстве $\{x \in \mathbb{R}^n: \|x- a\|_2 \leq \|x - b\|_2\}$. Функция ограничений имеет следующий вид:
    \begin{gather}\label{cons1_quasi}
        g(x) = \|x\|_2 + \max \{\langle -a, x\rangle, \|x\|_2\} - b, 
    \end{gather}
где $a \in \mathbb{R}^n$~--- фиксированный вектор и $b \in \mathbb{R}$, или имеет следующий вид:
\begin{gather}\label{cons2_quasi}
g(x) = \max_{i \in \overline{1, m}}\{g_i(x) = \langle \alpha_i, x \rangle - \beta_i \},  
\end{gather}
где $\alpha_i \in \mathbb{R}^n, \beta_i \in \mathbb{R}\, \forall i \in \overline{1, m}$. 

Для этого примера мы берем $a = \textbf{0} \in \mathbb{R}^n$, точка $b \in \mathbb{R}^n$ выбрана так, что $\|b-a\|_2 = 2$ и в качестве множества $Q$ выберем шар из $\mathbb{R}^n$ с центром в $a=\textbf{0} \in \mathbb{R}^n$ радиуса $1$.
\end{exl}

\begin{exl}[задача проектирования механических конструкций \cite{nesterov_2014}]\label{ex_lin_function}
Интересным примером приложений может быть использование схем с переключениями к задаче проектирования механических конструкций, которые могут быть сведены \cite{nesterov_2014} к оптимизационной задаче вида
\begin{gather}\label{prob_max}
    \max_{x \in Q} \langle \alpha, x \rangle, \quad \text{удовл.} \quad g(x) := \max_{i \in \overline{1, m}} \left\{ \pm\langle a_i, x\rangle - 1\right\} = \max_{i \in \overline{1, 2m}} g_{i} (x),
\end{gather}
где $g_i(x) = \pm \langle a_i, x\rangle - 1$ (у нас есть $2m$ функциональных ограничений $g_i(\cdot)$). 
\end{exl}

Выберем в качестве множества $Q$ для примеров~\ref{ex_quasi} и~\ref{ex_lin_function} евклидов шар из $\mathbb{R}^n$ с центром в $\textbf{0} \in \mathbb{R}^n$ радиуса $r$, a для примера~\ref{ex_geometric_prob} выбираем положительную часть этого шара.

Для примера~\ref{ex_geometric_prob} запускался алгоритм~\ref{alg:eps_sharp_min} (на практике, по крайней мере для рассмотренных примеров, заметной разницы в работе алгоритмов~\ref{alg:eps_sharp_min} и~\ref{alg:sharp_min_without_eps} нет, поэтому мы рассматриваем здесь только работу алгоритма~\ref{alg:eps_sharp_min}) со стартовой точкой $\left( \frac{1}{\sqrt{n}}, \ldots, \frac{1}{\sqrt{n}}\right) \in \mathbb{R}^n$ при $n=1000,\;m=100, r =1, \varepsilon = 10^{-3}, p = 5$. Коэффициенты в \eqref{cons_gemetric_prob} генерируются случайным образом с равномерным распределением на интервале $[0,1)$, и константы $b_i$  генерируются случайным образом с нормальным (гауссовским) распределением с математическим ожиданием в $0$, и среднеквадратичным отклонением, равным $1$. В этом случае оптимальное значение $f^* = 0$.

Для примера~\ref{ex_quasi} с ограничением вида \eqref{cons1_quasi}, запускались алгоритмы со стартовой точкой $x_0 = \left( \frac{1}{\sqrt{n}}, \ldots, \frac{1}{\sqrt{n}}\right) \in \mathbb{R}^n$ и $n = 10^5, r = 1$. Вектор $a$ и  константа $b$ в \eqref{cons1_quasi} сгенерированы случайным образом с равномерным распределением на интервале $[0,1)$. В случае ограничений вида \eqref{cons2_quasi} алгоритмы запускались с начальной точкой $\left( \frac{-1}{\sqrt{n}}, \ldots, \frac{-1}{\sqrt{n}}\right) \in \mathbb{R}^n$ при $n=1000,\;m=100, r = 1$. Коэффициенты $\alpha_i$ и $\beta_i$  в \eqref{cons2_quasi} генерируются случайным образом с равномерным распределением на интервале $[0,1)$. 

Для примера~\ref{ex_lin_function} мы запускали алгоритмы со стартовой точкой $x_0 = \left( \frac{1}{\sqrt{n}}, \ldots, \frac{1}{\sqrt{n}}\right) \in \mathbb{R}^n$ и $n = 1000, m = 100, r = 1, \varepsilon = 10^{-4}$. Вектор $\alpha$ в \eqref{prob_max} сгенерирован случайным образом с равномерным распределением на интервале $[0,1)$. Коэффициенты $a_i$ в \eqref{prob_max} генерируются случайным образом с нормальным (гауссовским) распределением с математическим ожиданием в $0$, и среднеквадратичным отклонением, равным $0.1$ и $1$.  

\begin{figure}[htp]
\centering
{\includegraphics[width=7.8cm]{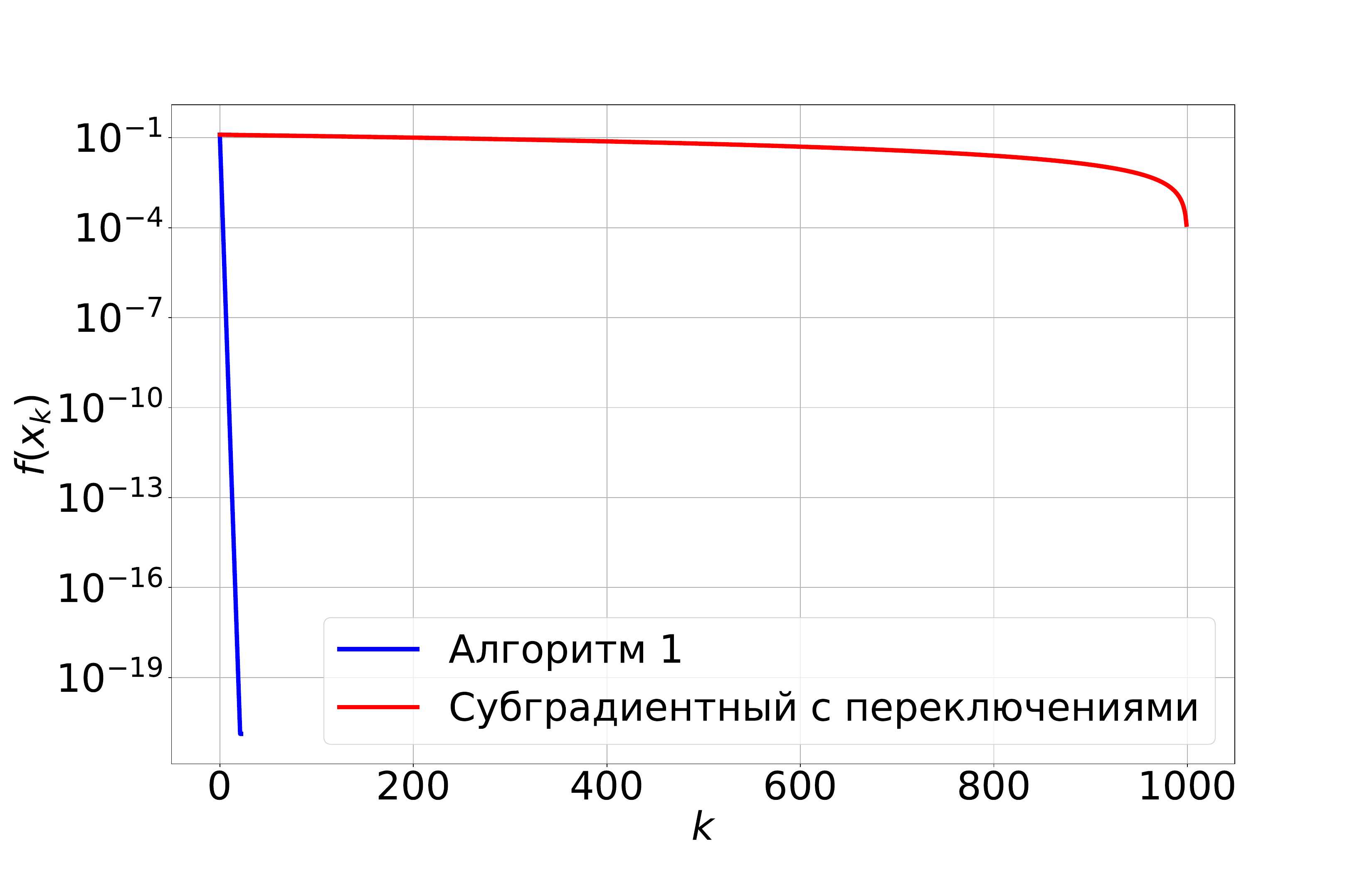} }
{\includegraphics[width=7.8cm]{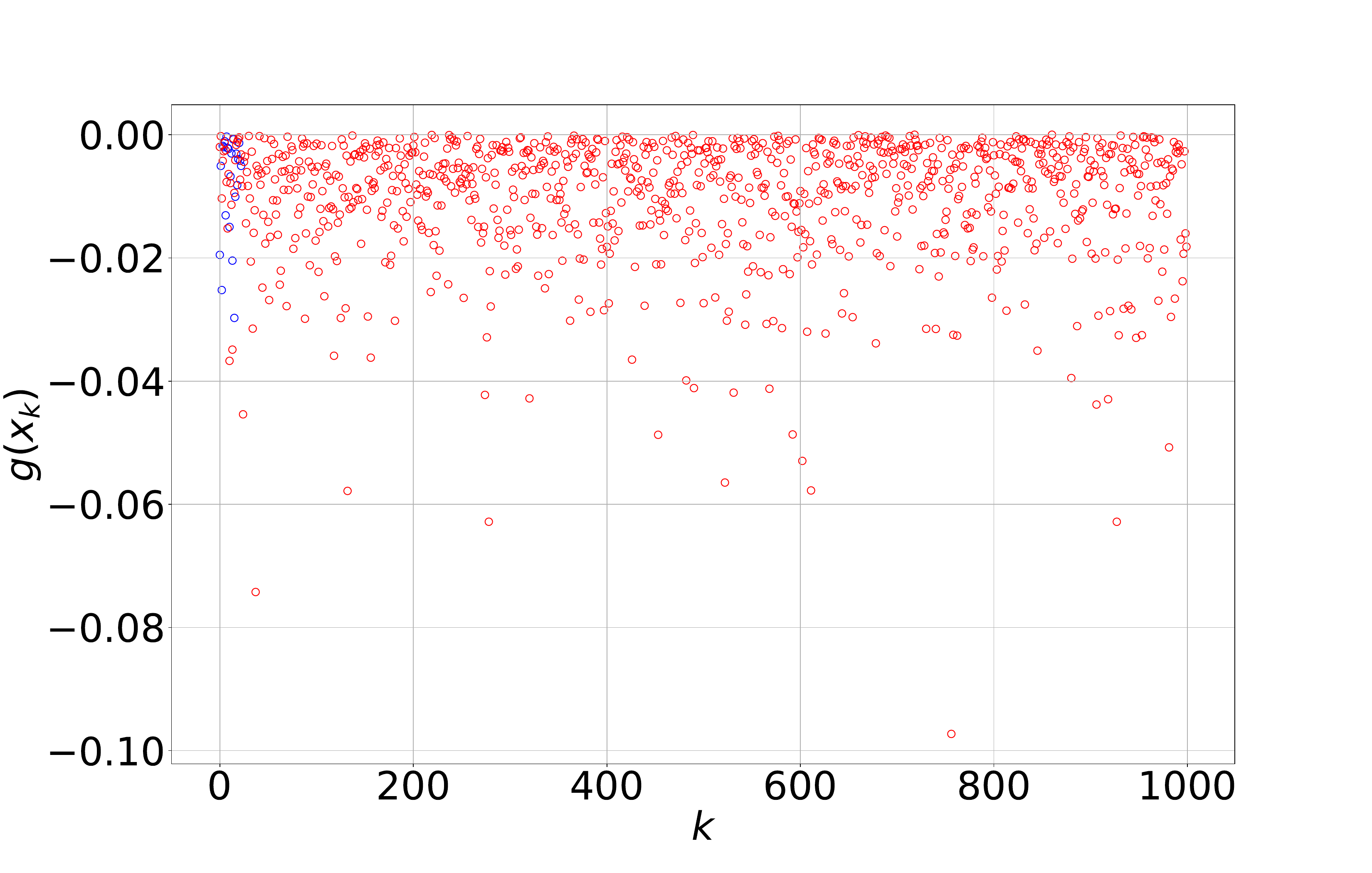} }
\caption{Результаты алгоритма~\ref{alg:eps_sharp_min} и субградиентного метода (алгоритм 3 из \cite{Ablaev2023}) с переключениями для примера~\ref{ex_geometric_prob}.  На этих рисунках показана динамика значений $f(x_k)$ \eqref{obj_lp} и значения $g(x_k)$ \eqref{cons_gemetric_prob}.} 
\label{res_geometric_problem}
\end{figure}

\begin{figure}[htp]
\centering
{\includegraphics[width=7.8cm]{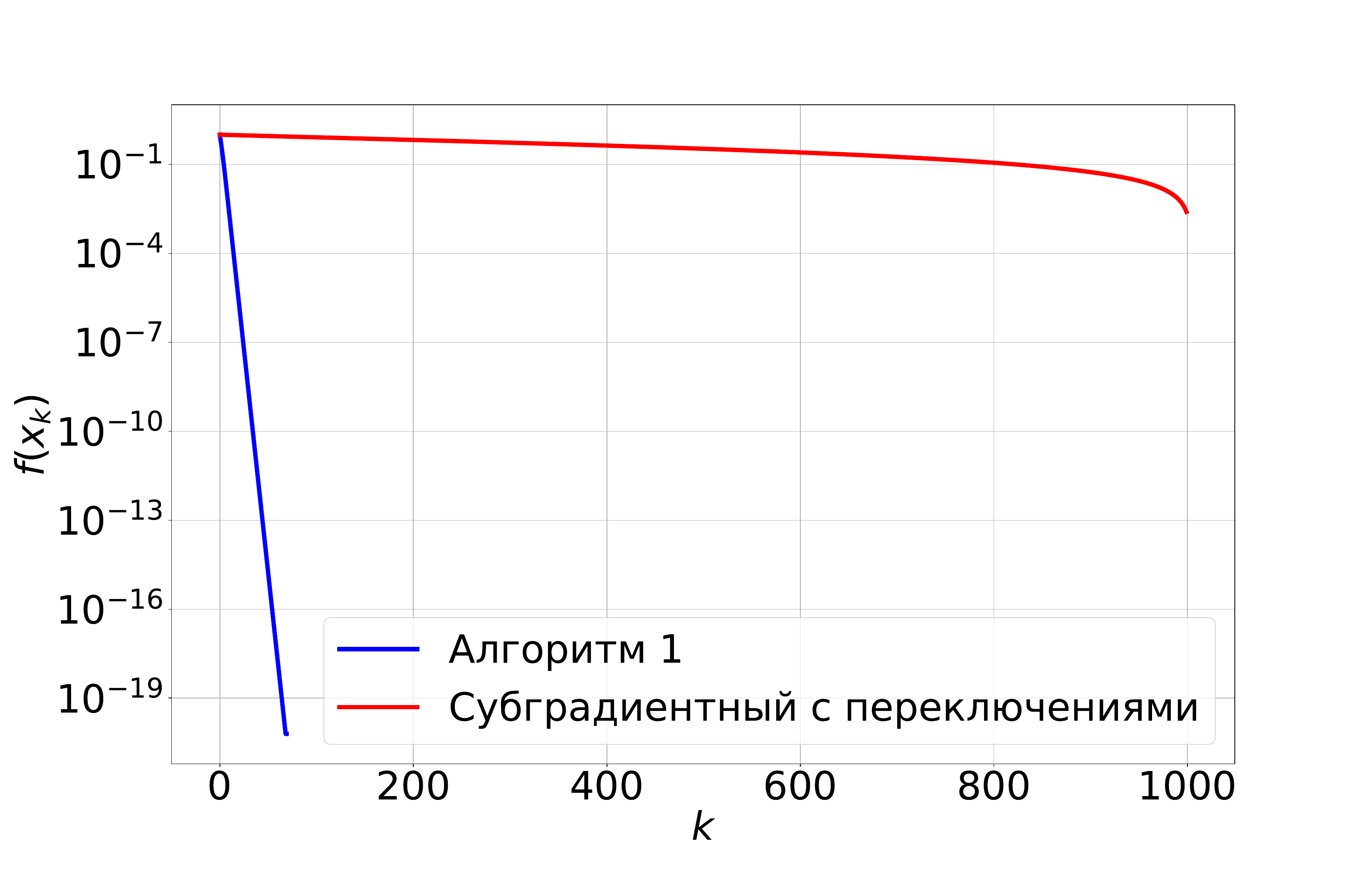} }
{\includegraphics[width=7.8cm]{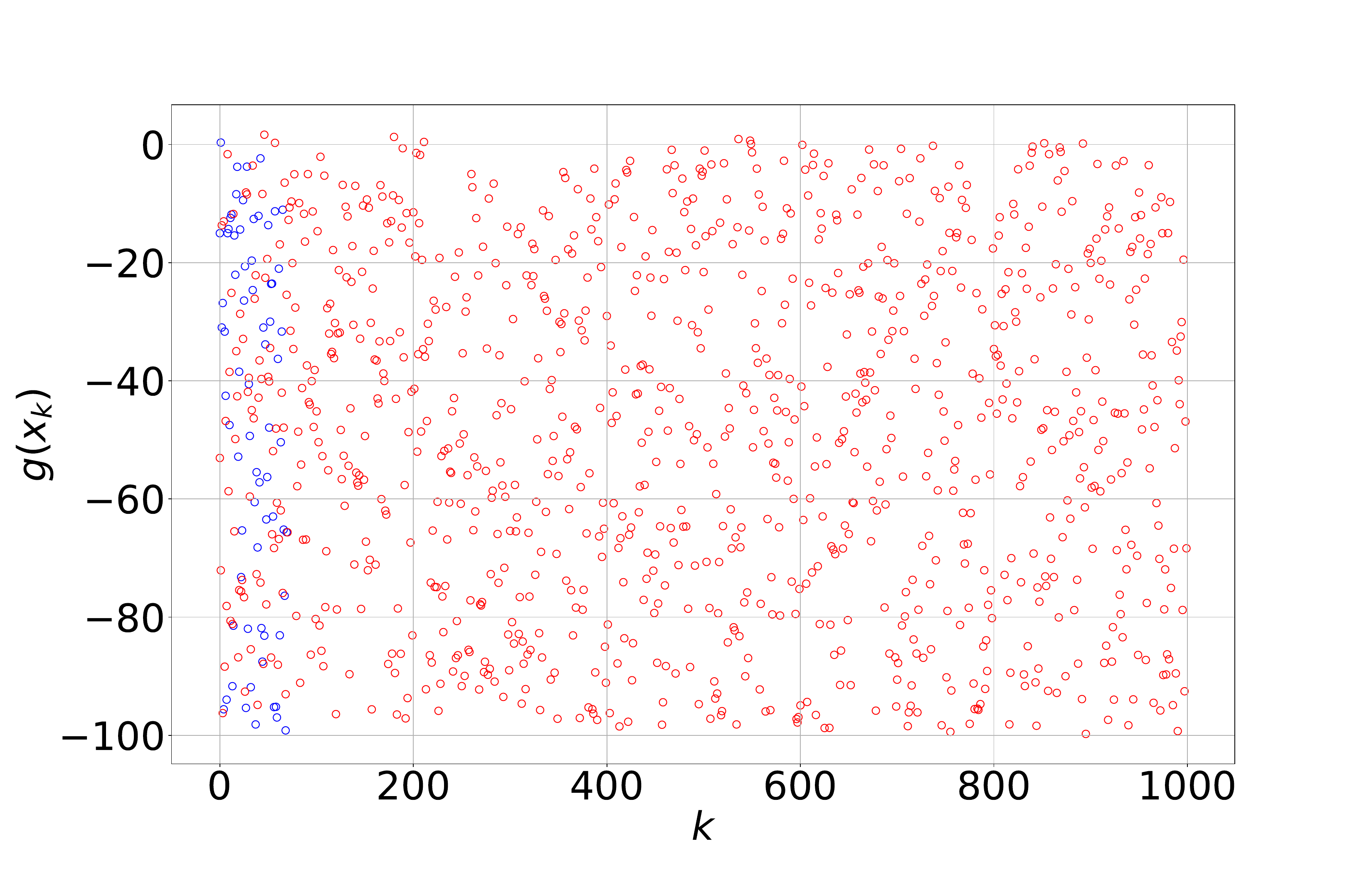} }
\caption{Результаты алгоритма~\ref{alg:eps_sharp_min} и субградиентного метода с переключениями (алгоритм 3 из \cite{Ablaev2023}) для примера~\ref{ex_quasi} с ограничениями \eqref{cons1_quasi}.  На этих рисунках показана динамика значений $f(x_k)$ \eqref{obj_quasi} и значения $g(x_k)$ \eqref{cons1_quasi}.}
\label{res_quasi_1}
\end{figure}

\begin{figure}[htp]
\centering
{\includegraphics[width=7.8cm]{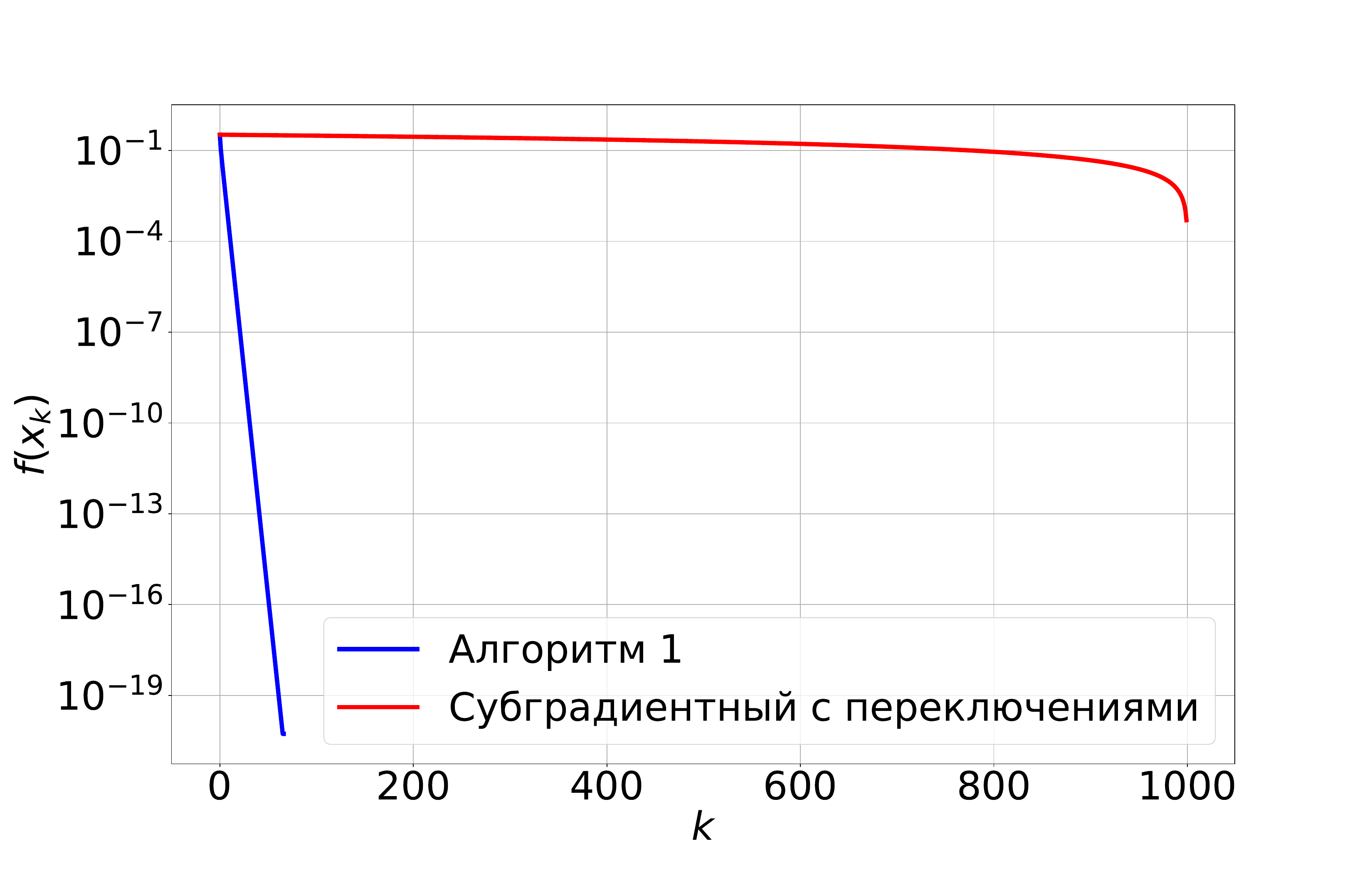} }
{\includegraphics[width=7.8cm]{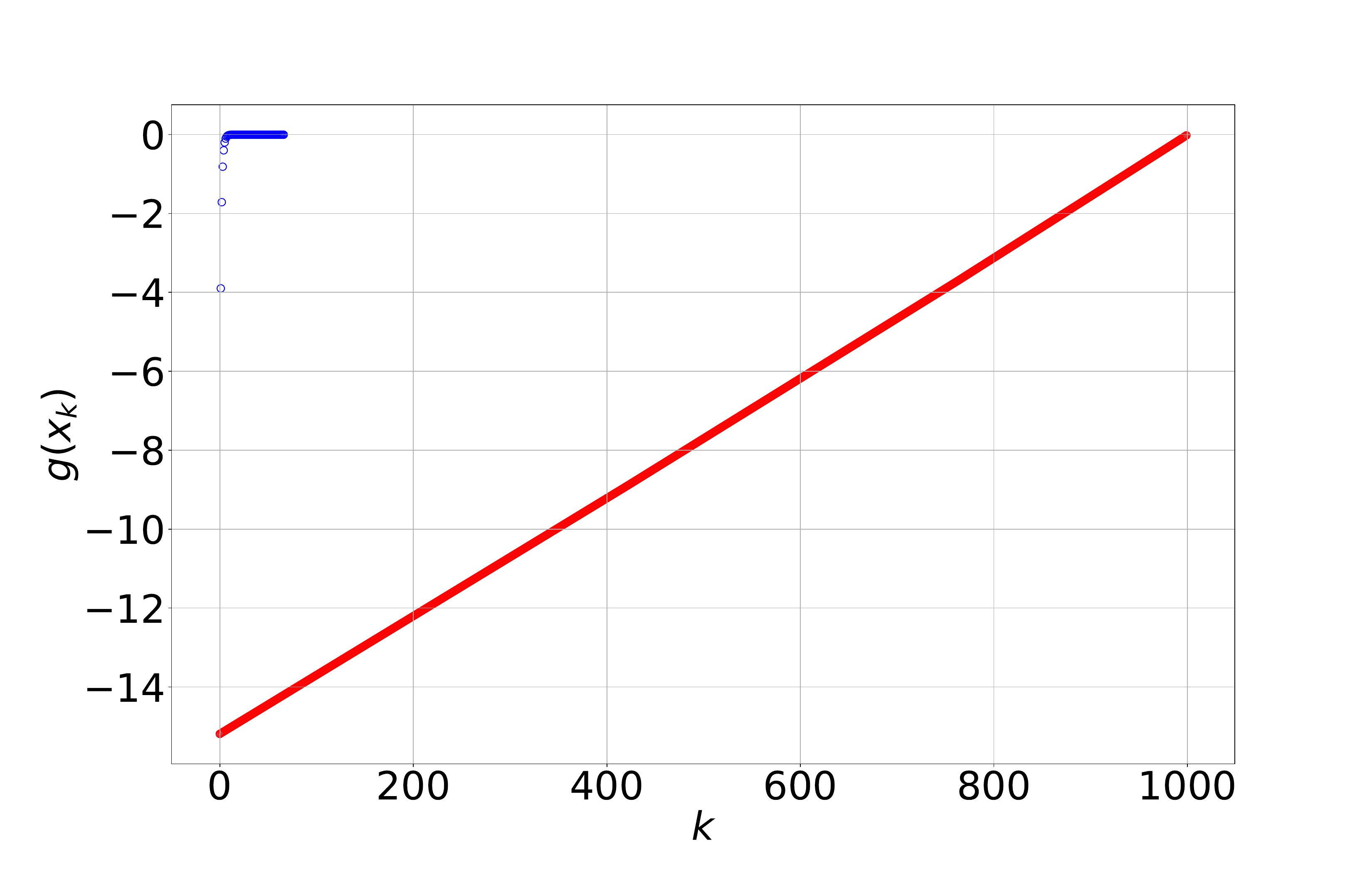} }
\caption{Результаты алгоритма~\ref{alg:eps_sharp_min} и субградиентного метода с переключениями (алгоритм 3 из \cite{Ablaev2023}) для примера~\ref{ex_quasi} с ограничениями \eqref{cons2_quasi}. На этих рисунках показана динамика значений $f(x_k)$ \eqref{obj_quasi} и значения $g(x_k)$ \eqref{cons2_quasi}.}
\label{res_quasi_2}
\end{figure}

\begin{figure}[htp]
\centering
{\includegraphics[width=7.8cm]{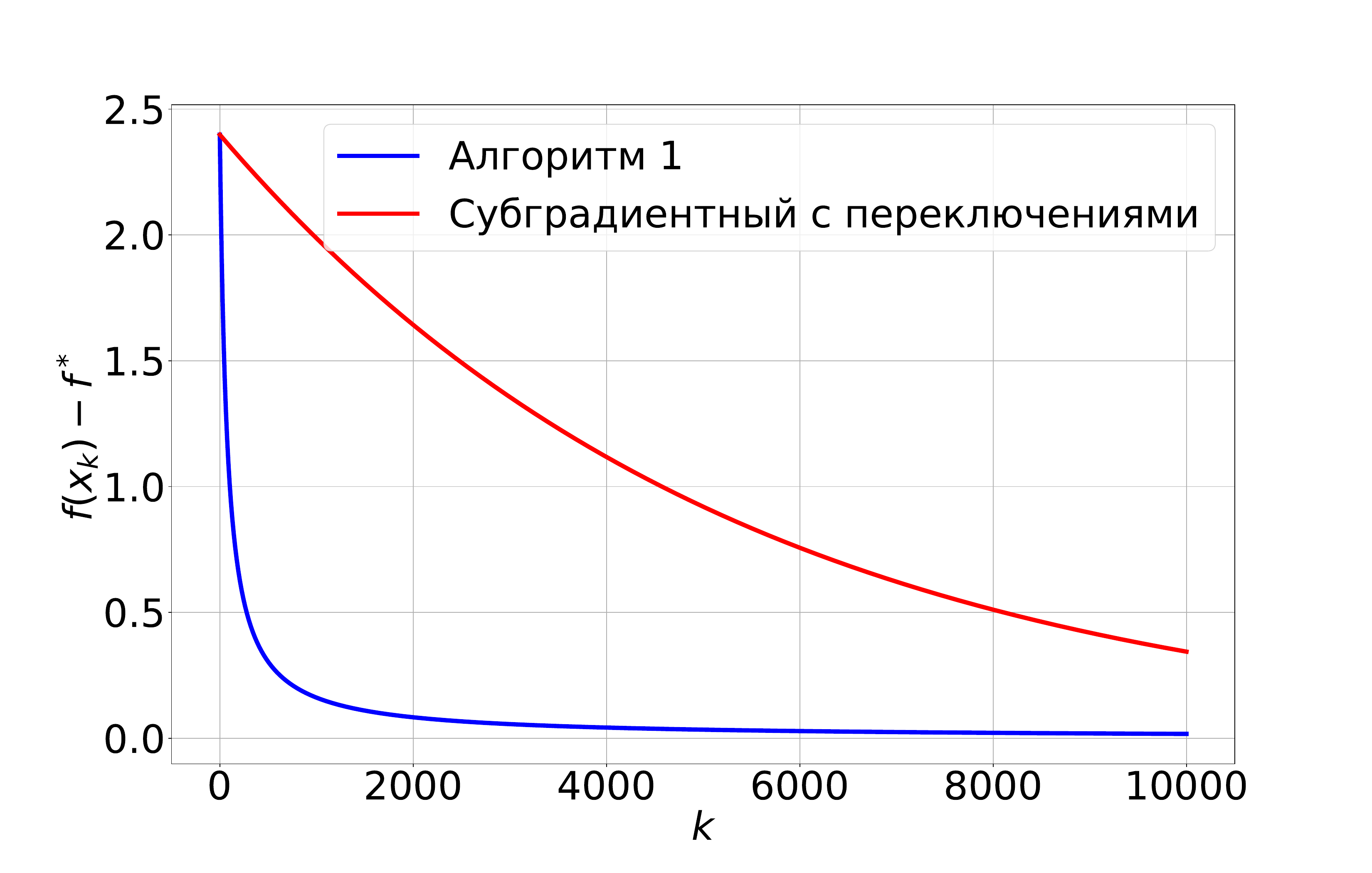} }
{\includegraphics[width=7.8cm]{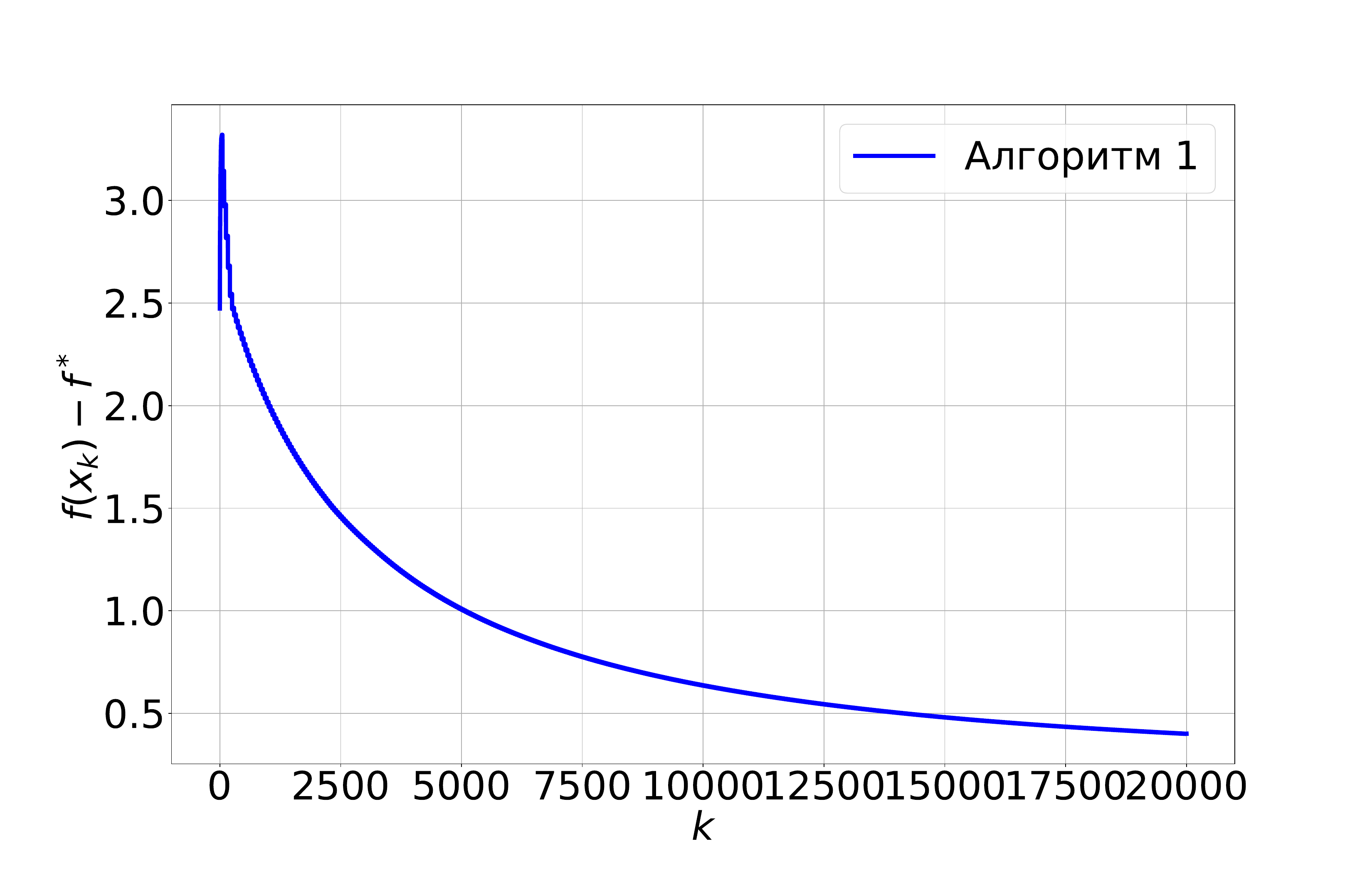} }
\caption{Результаты алгоритма~\ref{alg:eps_sharp_min} и субградиентного метода с переключениями (алгоритм 3 из \cite{Ablaev2023}) для примера~\ref{ex_lin_function}  (задача \eqref{prob_max}). На этих рисунках показана динамика значений $f(x_k) - f^*$, когда коэффициенты $a_i$ в \eqref{prob_max} генерируются случайным образом с нормальным распределением с центром, равным $0$, и стандартным отклонением, равным $0.1$, (слева) и с центром, равным $0$, и стандартным отклонением, равным $1$, (справа).}
\label{res_max_problem}
\end{figure}

Результаты сравниваемых алгоритмов (алгоритм~\ref{alg:eps_sharp_min} и алгоритм 3 из \cite{Ablaev2023} <<субградиентный метод с переключениями>>) представлены на рис.~\ref{res_geometric_problem},~\ref{res_quasi_1},  ~\ref{res_quasi_2} и~\ref{res_max_problem}. 

Эти результаты отражают значения целевой функции в примерах~\ref{ex_geometric_prob} и~\ref{ex_quasi}, а также значения функции ограничений в точках на каждой итерации. А результаты на рис.~\ref{res_max_problem} (для примера~\ref{ex_lin_function})  отражают значении $f(x_k) - f^*$, когда коэффициенты $a_i$ в \eqref{prob_max} генерируются случайным образом с нормальным распределением с математическим ожиданием $0$ и стандартным средним квадратичным отклонением, равным $0.1$, (слева) и с математическим ожиданием $0$ и стандартным средним квадратичным отклонением $1$ (справа).

Из рис.~\ref{res_geometric_problem},~\ref{res_quasi_1},~\ref{res_quasi_2} и ~\ref{res_max_problem} мы видим эффективность предложенного алгоритма~\ref{alg:eps_sharp_min}, с помощью этого алгоритма мы можем достичь решения задачи с очень высокой скоростью по сравнению с другими схемами, такими как алгоритм 3 из \cite{Ablaev2023}, которые медленно сходятся к решению.  Также на рис.~\ref{res_max_problem} (справа) мы видим результаты только для алгоритма~\ref{alg:eps_sharp_min}, потому что для случая, когда коэффициенты $a_i$ в \eqref{prob_max} генерируются случайным образом с нормальным распределением с математическим ожиданием $0$ и срденим квадратичным отклонением $1$, по субградиентному алгоритму с переключением до $20000$ итераций нет ни одного продуктивного шага.  


Сравним теперь скорость сходимости метода для шага вида $h_k^f = \frac{\epsilon}{\|\nabla f(x_k)\|_2^2}$ и шага Поляка из алгоритма~\ref{alg:eps_sharp_min}: $h_k^f = \frac{f(x_k) - \overline{f}}{M_f \|\nabla f(x_k)\|_2}$. Шаг для ограничений же у обоих методов поставим одинаковый: $h_k^g = \frac{1}{\|\nabla g(x_k)\|_2}$. Далее, рассмотрим следующий пример.

\begin{figure}[htp]
\centering
{\includegraphics[width=7.8cm]{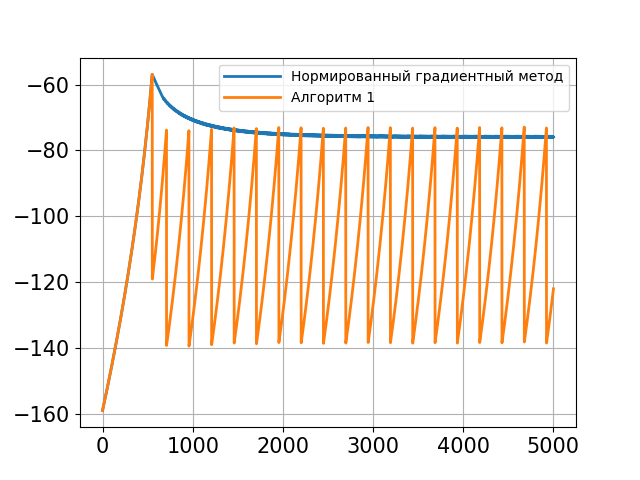} }
\caption{Результаты сходимости алгоритма~\ref{alg:eps_sharp_min} с шагом Поляка и с нормированным градиентным шагом для примера \eqref{ex_polyak_and_norm_grad}.}
\label{res_ex_polyak_and_norm_grad}
\end{figure}

\begin{exl}\label{ex_polyak_and_norm_grad}
Предположим, что $x \in \mathbb{R}^n_{++}$, и определим следующую целевую функцию: 
    \begin{gather}\label{ex_polyak_and_norm_grad_cnstr_f}
        f(x) = - \sqrt{x^\top a}, 
    \end{gather}
    где $a \in \mathbb{R}^n_{++}$~--- фиксированный вектор, а функция ограничений имеет следующий вид:
    \begin{gather}\label{ex_polyak_and_norm_grad_cnstr}
        g(x) = \sum_{i=0}^n x \log{ (x_i / a_i) } - x_i - a_i, 
    \end{gather}
т.~е. $g(x)$ отражает расстояние Кульбака--Лейблера между $a$ и $x$. Положим $g(x) \leq 1000$ и размерность задачи $n=10^5$.

Для приближённого решения задачи минимизации $f(x)$ с ограничениями $g(x)$ на множестве $\mathbb{R}^n_{++}$ рассчитаем решение в библиотеке cvxpy. Воспользуемся этим примерным решением для расчёта шага Поляка на продуктивных итерациях. Из рис.~\ref{res_ex_polyak_and_norm_grad} видим, что метод с шагом Поляка сильно превосходит шаг вида $h_k^f = \frac{\epsilon}{\|\nabla f(x_k)\|_2^2}$ по скорости сходимости, его убывание гораздо резче, но затем, после достаточного убывания, оба шага отбрасывает назад непродуктивными шагами.

\end{exl}

\paragraph{6. Заключение}

В статье исследованы некоторые субградиентные методы с переключениями по продуктивным и непродуктивным шагам для задач минимизации квазивыпуклых липшицевых функций с ограничениями-неравенствами двух предположений об остром минимуме \cite{Stonyakin2023,Lin2020}. Для такого класса задач развивается подход \cite{Ablaev2023} на базе субградиентных методов с переключениями по продуктивным и непродуктивным шагам. Однако важно, что вместо рестартов по параметру острого минимума используется процедура регулировка шага в такого типа сyбградиентных методах, аналогичная подходу Б.Т. Поляка \cite{Polyak1969}. Известно, что использование в методах градиентного типа шага Б.\,Т.~Поляка довольно популярно (см., например, \cite{Hazan2019, Loizou2021}). Совсем недавно предложены и изучены различные современные варианты шага Поляка \cite{Wang2023, Devanathan2023, Abdukhakimov2023}. При этом однако естественно возникает проблема доступности информации о минимальном значении целевой функции, что ещё более важно для задач с дополнительными ограничениями. В статье проработан этот вопрос и предложен подход к такого типа процедуре регулировки шага при условию с неточной информацией об оптимальном значении целевой функции задачи. Получены условия, про которых для такого типа шагов можно ожидать сходимость метода со скоростью геометрической прогрессии в окрестность множества точных решений. Описана методика повышения скорости работы данного алгоритма в случае задач со многими ограничениями-неравенствами, основанная на отказе от прохода по всем ограничениям.Доказаны оценки качества выдаваемого методом решения для предложенных алгоритмов в зависимости от точности информации о минимальном значении $f^*$. Выполненные вычислительные эксперименты, особенно для задач геометрического программирования и проектирования механических конструкций показали хорошую эффективность предложенного метода в решении прикладных задач. Существенный момент в том, что использование предлагаемой в настоящей работе процедуры регулировки шага позволяет при реализации уйти от необходимости знания параметра острого минимума, который на практике часто довольно проблематично оценить. Вообще такого типа шаги потенциально можно применять к любым оптимизационным задачам с ограничениями-неравенствами. Поэтому представляется, что исследованные в статье методы перспективны для использования на классах самых разных минимизационных задач с ограничениями.


\end{document}